\documentclass[11pt]{article}
\usepackage{multirow}
\usepackage{amsfonts, amssymb, graphicx, a4}
\usepackage{amsmath,dsfont}
 \usepackage{color}
 \usepackage{hyperref}

\newtheorem{theorem}{Theorem}[section]
\newtheorem{lemma}[theorem]{Lemma}
\newtheorem{proposition}[theorem]{Proposition}

\newtheorem{definition}[theorem]{Definition}
\newtheorem{corollary}[theorem]{Corollary}
\newtheorem{remark}[theorem]{Remark}
\newtheorem{remarks}[theorem]{Remarks}

 \newenvironment{proof}{\noindent {\bf Proof.\,}
 }{\hspace*{\fill}$\qed$\medskip}

\def \P {{\mathbb P}}
\def \E {{\mathbb E}}

\newcommand{\cX}{{\cal X}}

\def \a {{\alpha}}

\def \e {{\varepsilon}}

\def \r {{\rho}}
\def \l {{\lambda}}
\def \m {{\mu}}

\def \s {{\sigma}}

\def \g {{\gamma}}
\def \t {{\tau}}

\def \d {{\delta}}
\def \p {{\pi}}

\def\eqref#1{(\ref{#1})}

\newcommand{\be}[1]{\begin{equation}\label{#1}}
\newcommand{\ee}{\end{equation}}

\newcommand{\bl}[1]{\begin{lemma}\label{#1}}
\newcommand{\el}{\end{lemma}}

\newcommand{\br}[1]{\begin{remark}\label{#1}}
\newcommand{\er}{\end{remark}}

\newcommand{\bt}[1]{\begin{theorem}\label{#1}}
\newcommand{\et}{\end{theorem}}

\newcommand{\bd}[1]{\begin{definition}\label{#1}}
\newcommand{\ed}{\end{definition}}

\newcommand{\bcl}[1]{\begin{claim}\label{#1}}
\newcommand{\ecl}{\end{claim}}

\newcommand{\bp}[1]{\begin{proposition}\label{#1}}
\newcommand{\ep}{\end{proposition}}

\newcommand{\bc}[1]{\begin{corollary}\label{#1}}
\newcommand{\ec}{\end{corollary}}

\newcommand{\bi}{\begin{itemize}}
\newcommand{\ei}{\end{itemize}}

\newcommand{\ben}{\begin{enumerate}}
\newcommand{\een}{\end{enumerate}}
 \newcommand{\bpr}{\begin{proof}}
 \newcommand{\epr}{\end{proof}}

\def \qed {{\square\hfill}}

\def \P {{\mathbb P}}
\def \E {{\mathbb E}}

\def\1{\rlap{\mbox{\small\rm 1}}\kern.15em 1}

\def\build#1_#2^#3{\mathrel{\mathop{\kern 0pt#1}\limits_{#2}^{#3}}}
\def\tend#1#2#3{\build\hbox to 12mm{\rightarrowfill}_{#1\rightarrow
#2}^{#3}}

{}
\def\converge#1#2#3{\build\hbox to
15mm{\rightarrowfill}_{\hbox{\scriptsize #3}}^{#1\rightarrow #2}}
\def\converg#1#2#3{\build\hbox to
15mm{\rightarrowfill}_{\hbox{\scriptsize #3}}^{#1\uparrow #2}}

\newcommand{\eproof}{\hspace*{\fill}$\qed$}

\begin{document}

\title{Strong  times and first hitting. }

\author{F.\ Manzo,
\footnote{ Supported by Dipartimento di Matematica e Fisica,
 Universit\`a di Roma Tre
 }\\
E.\ Scoppola
 \footnote{
 Dipartimento di Matematica e Fisica,
 Universit\`a di Roma Tre,
 Largo S.\ Leonardo Murialdo 1,
 00146 Rome, Italy
 }
}
\maketitle
\begin{abstract}

We generalize the notion of  strong stationary time and we give a representation
formula for the hitting time to a target set in the general case of  non-reversible Markov processes.
\end{abstract}


\section{Introduction }

This work was originally motivated by the study of the first hitting
to rare sets for ergodic Markov chains. 
Our aim was to provide a unifying
language for different approaches to the problem, focusing on the
link between rarity and exponentiality,  in particular for metastable systems.

Fundamental results about the connection between rarity and exponentiality
trace back to the late '70s, in particular to the book of Keilson
\cite{Kei79}.  
The ideas and methods used (reversibility, complete-monotonicity,
spectral decomposition, fundamental matrix, potential theory) were developed
and generalized in a subsequent series of papers, among which we want
to mention \cite{A,AB1,AB2,B99,FL}.
For a discussion of early literature, we refer to \cite{AG00}.

Central results in this approach, concern the distance between $e^{-t}$ and the tail distribution 
of the ratio between the hitting time $\tau_G$ and its mean starting from 
the invariant measure $\pi$
\be{th12}  \P_\pi
  \left( \frac{\tau_G}{\E_\pi \tau_G} > t \right).
\ee

Explicit, time-uniform, bounds of this distance were given, starting from \cite{A}, in terms of the ratio between a ``local relaxation time" $R$ and
$\E_\pi(\tau_G)$.
In the literature, the role of this local relaxation time $R$  was played by many different times such as 
the mixing time
\cite{A},  the relaxation time (i.e., the inverse
of the smallest non-zero eigenvalue of the generator of the dynamics) in \cite{AB1}, the inverse of the spectral radius of the sub-markovian generator restricted outside $G$ \cite{BG} and some other, sometimes model-dependent, choices.
Heuristically, when $R \ll \mathbb{E}_{\pi}\tau_{G}$,
the system relaxes to a local equilibrium before attempting to reach $G$.
Extensions to other initial measures 
(in particular to the conditional equilibrium measure $\pi(\cdot|\cX\backslash G)$ and to the quasi-stationary measure \cite{AB1}), or to non-reversible settings were given over the years.   

Rather independently, in the early eighties the statistical mechanics community began
to study metastability 
as a dynamical phenomenon,
in the framework of discrete models with Glauber
dynamics. 
The regimes studied in this context were not as general as
that in the above mentioned papers but, on the other hand, the physical meaning
of the exponential behavior was very transparent.
In this case the target set $G$ is often the ``basin of attraction" of the stable state and
the role of the relaxation time $\E_\p\t_G$ is taken by $\E_\m\t_G$ where $\m$ is some ``metastable 
measure" concentrated ouside this basin.
We refer to \cite{AG00} and \cite{FMNS} for a discussion 
and comparison of the different approaches to prove exponential behavior in metastable system.
We emphasize the fact that a 
characterization of ``good metastable" initial measures, that give rise to the exponential behavior of the relaxation time, is of primary interest in this series of works.

Let us give a very quick heuristic. In the easiest asymptotic regime, the system is trapped into an energy well and spends in the bottom of this well most of the time before
reaching the boundary and exit. 
If this bottom is represented by a single point,
each return time the process looses
memory of its past. 
This renewal property gives rise to the exponential
behavior. 
Large-deviation methods, renormalization ideas, coupling, simulated
annealing techniques and, more recently, potential theory and martingales have been used
to extend this picture to physically more interesting settings,
provided the renewal properties of a metastable
point are strong
enough.

Here we generalize this language, by
``changing the renewal point for a renewal measure''. 
The obvious candidate for this role is the
\emph{quasi-stationary measure} (see \cite{DS}), namely: 
\begin{equation}
\mu^{*}(\cdot):=\lim_{t\to\infty}\mathbb{P}\left(X_{t}=\cdot \ | \ t<\tau_{G}\right).\label{mu*}
\end{equation}
 Indeed, it is easy to see that the evolution starting from this measure
is exponential, in the sense that 
\begin{equation}
\P_{\mu^*}\left( X_t = y\right) =:
\mu_{t}^{\mu^{*}}(y)=\begin{cases}
\lambda^{t}\mu^{*}(y) & \text{ if }y\not\in G\\
(1-\lambda^{t})\omega(y) & \text{ if }y \in G
\end{cases}\label{mumu*}
\end{equation}
where $\lambda$ is the largest eigenvalue of the sub-markovian matrix 
 obtained from $P$ by canceling out the entries in $G$,
and where $\omega(y):=\P(\tau^{\mu^*}_G = \tau_y)$
is the probability that $y$ is the hitting point.

We will refer to the measure $\mu_{t}^{\mu^{*}}$ as the \emph{squeazing quasi-stationary measure}.

The idea is then to control how close is $\mu_{t}^{\alpha}(y):=\P\left( X^{\alpha}_t = y\right) $ to
$\mu_{t}^{\mu^{*}}(y)$ when $\alpha$ is some other starting measure.

In order to proceed in this direction, we introduce a sort of ``hitting
time to a measure'' by generalizing the notion of \emph{strong stationary
time}, introduced in \cite{AD1} (under the name of \emph{strong uniform
time}).

Using the strong time language, we are able to give a representation formula for the 
probability 
$\P\left( \tau^\a_G > t \ ; \ X^\a_t=y \right)$
in terms of events concerning these strong times.
This representation formula gives a probabilistic interpretation of the errors in the exponential approximation and is very explicit about the role of the initial measure $\alpha$.

Let us recall, from \cite{AD1}, the following 
\begin{definition}

A randomized stopping time $\tau^\alpha_\pi$ is a \emph{Strong Stationary Time } (SST in the following) for the
Markov chain $X^\a_t$ with starting distribution $\a$ and stationary measure $\p$,  if
\[
\mathbb{P}\left(X^{\alpha}_t=y,\tau^{\alpha}_\pi=t\right)=\pi(y)\mathbb{P}\left(\tau^{\alpha}_\pi=t\right).
\]
\end{definition}
SSTs were introduced in  \cite{AD1}, where their existence was proved.  
The proof also shows that the fastest SST  is distributed according to the separation between the measure
at time $t$ and the stationary measure.

Explicit constructions of SSTs can be done in very particular cases, in  one dimension  or in  very symmetric systems (\cite{AD2,DF90}), where these constructions were used e.g. to show cutoff behavior.

Separation itself is not the easiest notion of distance between measures to compute, but the separation between the measure
at time $t$ and the stationary measure has the remarkable property of being submultiplicative and makes it usable to give exponential bounds.

In our point of view, strong times provide a new language to describe the approach to equilibrium or, in our case, to quasi-stationarity.
Actually, inspired by metastability, we will consider the hitting time $\tau_G$ as the decay time
of the metastable state. It is then natural to assume $G$ an absorbing state and ergodicity on $\cX\backslash G$.
We are not interested in finding explicit constructions:
all we want to do is  to use exponential bounds on the separation in order to estimate the tail distribution of the hitting time.

Let us mention that the idea of using a strong time that somehow
catches the arrival to the quasi-stationary measure is not new in
the literature; in \cite{DM09}, for a birth-and-death process starting from $0$, in a particular
regime, the authors construct what they call a \emph{``strong
quasi-stationary time''} for this purpose. 
\bigskip

Our approach is different under two fundamental aspects:
\begin{enumerate}
\item Our notion of Conditionally Strong Quasi Stationary Time is completely
general and its existence does not require any additional assumptions
besides ergodicity of the stochastic matrix outside $G$.
The prize to pay is that, in general,
 we cannot construct explicitly such times.
\item Our target is not a fixed measure but a family of measures indexed
by the time $t$. The reason behind our choice is that while the measure
$\mu^{*}$ is concentrated outside $G$, in general the evolved measure
$\mu_{t}^{\alpha}$ is not. 
On the contrary, when $G$ is
absorbing, $\mu_{t}^{\alpha}$ concentrates on $G$. Therefore, for
general models and general starting states, there is no hope to reach
$\mu^{*}$ at a positive time. 
\end{enumerate}
A natural candidate for the role of
``target measure'' is instead the \emph{``squeazing measure''}
$\mu_{t}^{\mu^{*}}$ or, more in general,
a family of measures  $\m_t$, with the property
$\m_{t+1}(x)=\sum_{y \in \cX}\m_t(y)P_{x,y}$. 

This choice allows to define properly a strong
time $\tau_{\m}^{\alpha}$ such that 
\be{SQST}
\mathbb{P}\left(X^{{\alpha}}=y,\tau_{\m}^{\alpha}=t\right)=\m_{t}(y)\mathbb{P}\left(\tau_{\m}^{\alpha}=t\right).
\ee
Unfortunately, as we will see, this time decays in a time of order $\E(\tau^\a_G)$
and it is too large for the applications we have in mind. 
The reason is, $\tau_{\m}^{\alpha}$
gives the same role to the points in $G$ and outside $G$. A good
``local relaxation time'' instead, should regard only what happens
outside $G$.

For this reason, it is natural to consider a conditional time:
\begin{definition}
A randomized stopping time $\tau^\alpha_*$ is a \emph{conditionally-strong quasi-stationary time}  (CSQST in the following) 
if for any $y \not \in G$, 
\be{CSQST}
\P \left( X^\alpha_t = y, \tau^\alpha_*=t \ | \ t < \tau^\alpha_G \right) 
= 
\mu^*(y) \P \left(  \tau^\alpha_*=t \ | \ t < \tau^\alpha_G \right) 
\ee
\end{definition}
or, in other words, 
\be{defCSQST'}
 \P(X^{\alpha}_t=y,\, \t^{\alpha}_*=t)=\m^*(y)\P(\t^{\alpha}_*=t<\t^{\alpha}_G).
 \ee
The idea is  to use this CSQST in the decomposition
\be{basic}
 \P(\t^{\alpha}_G>t)=
 \P(\t^{\alpha}_G>t \ ; \ \t^{\alpha}_{*}\le t)
 +\P( \t^{\alpha}_{*,G}> t),
\ee
 where $\t^{\alpha}_{*,G}={\t^{\alpha}_G\wedge \t^{\alpha}_*}$. 

The event in the first term in the r.h.s. of (\ref{basic}) can be read as the probability that the process reaches $G$ after reaching the ``metastable equilibrium". 
Since in our setting, the metastable equilibrium is related to the quasi-stationary measure $\mu^*$, we easily get exponential bounds.
 
Its counterpart is the event that the process stays away from $G$ without reaching the ``metastable equilibrium". 
From (\ref{basic}), we obtain a probabilistic interpretation of the error term in the exponential bound.
 The role of the local relaxation time, in our approach,  is being played by $\t^{\alpha}_{*,G}$.

Both terms in r.h.s of (\ref{basic}) have exponential decay for large $t$. 
The exponential behavior of  $\P(\t^{\alpha}_G>t)$ emerges when the last term decays faster than the other and can be neglected.

The role of the initial measure $\alpha$  can be further clarified by the introduction of a ``time-shift":
different starting measures can help or hinder achieving $G$.  
Asymptotically, this fact results in a time-shift  $\delta_\alpha$, and our choice of CSQST corresponds to $\rho_t = \mu^{\mu^*}_{t+\delta_\alpha}$  
in (\ref{SQST}).

Indeed, just like {\it fastest } SSTs are related to the separation between the measure at a given time and the stationary measure, 
{\it ``minimal"} CSQSTs are related to the separation
\be{stilde}
  \tilde{s}^\alpha(t):=  \max_{y \not \in G} \  1-\frac{\mu^{\alpha}_{t}(y)}{\mu^{\mu^*}_{t+\delta_\alpha}(y)},
\ee
which quantifies a sort of distance between $\mu^\alpha_t$ and $\mu^{\mu^*}_{t+\delta_\alpha}$ by taking care only of the points outside $G$.

In subsection \ref{localch}
 we will show that $\tilde s^\alpha$
is the separation between the evolution of an auxiliary Markov chain in $\cX\backslash G$ and its stationary measure.  
Therefore $\max_{\alpha} \tilde s^\alpha$ it is submultiplicative and decays exponentially in time.

By using $\tilde s^\alpha$, (\ref{basic}) can be rephrased in a more usable representation formula as:
\be{eqmaint}
 \P(\t^{\alpha}_G>t)=\l^{t+\d_{\alpha}}
(1-\tilde s^\alpha (t))
 +\P( \t^{\alpha}_{*,G}> t).
\ee
We refer to subsection \ref{mainres}
for more precise statements.
\bigskip

{\bf Outline of the paper.} In subsection \ref{1.1} the general setting and definitions are fixed. In subsection \ref{localch}
we introduce a local chain on $\cX\backslash G$ which will be crucial in our discussion, while in subsection
\ref{mainres} our main results are stated. Section \ref{pf_main0} 
is devoted to the introduction of the central object in this paper: the generalization of strong stationary times to other target evolving measures. 
Subsection \ref{s1}
contains the proof of Theorem \ref{main0} and  subsection \ref{s2} the construction
of these strong times with an auxiliary chain. Section \ref{s_main1} contains the proof of Theorem \ref{main1}
and Section \ref{s_main2} the proof of Theorem \ref{main2}.
Concluding remarks and future perspectives are discussed in Section \ref{s5}.

\subsection{General setting and definitions }
\label{1.1}
We collect in this subsection definitions and notations used in the paper.
\begin{itemize}
\item
{\bf Process:}
we will consider  discrete time Markov chains  $\{X_t\}_{t\in \mathbb N}$  on a countable state space $\cX$.
We  denote  by $P(x,y)$ the
transition matrix
  and  by $\m^x_t(\cdot)$ the  measure at time $t$, starting at $x$, i.e.,  
$\m^x_t(y)\equiv\P(X^x_t=y)=P^t(x,y)$, for any $y\in \cX$.  More generally given an initial distribution
${\a}$ on $\cX$
$$
\m^{\a}_t(y)=\P(X^{\a}_t=y)=\sum_{x\in\cX}{\a}(x)P^t(x,y)
$$
Starting conditions (starting state $x$ or starting measure $\a$) will be denoted by a  
superscript in random variables  (i.e., $X^x_t$, $X^\a_t$,  $\t^x$, $\t^\a$...).

Let $G\subset \mathcal{X}$ be a target set and  $\tau_{G}$   its first hitting time 
$$\tau_{G}:=\min\{t\ge0\:;\: X_{t}\in G\}.$$
We will study   the process  $\{X_t\}_{t\in \mathbb N}$ up to time  $\tau_{G}$,  so it is not restrictive to assume that states in $G$ are absorbing.
Let $A:=\mathcal{X}\backslash G$, we assume ergodicity on $A$. More precisely,
 denoting by $[P]_A$  the sub-stochastic matrix obtained by $P$ by restriction to $A$
$$
[P]_A(x,y)=P(x,y)\ge 0 \qquad \forall x,y\in A,\qquad \sum_{y\in A}[P]_A(x,y)\le 1,
$$
we suppose $[P]_A$ a {\it primitive matrix}, i.e., there exists an integer $n$ such
that $\big([P]_A\big)^n$ has strictly positive entries.
\item
{\bf  Quasi-stationary measure on $A$:}
by the Perron-Frobenius theorem  it
 can be proved that there exists $\l<1$ such that $\l$ is the spectral radius of
  $[P]_A$ and there exists a unique non negative left eigenvector of $[P]_A$ corresponding to $\l$, i.e.,
  \be{eqmqs}
  \m^*[P]_A=\l \m^*.
  \ee
 so that  we get immediately 
 $$\P\left(\tau^{\mu^*}_{G}>t\right)=\l^t.$$
 \item
 {\bf Evolving measures:}
we call  {\it  evolving measure }
  any  family of measures  $\{\m_t\}_{t\in \mathbb N}$,  on $\cX$, 
  such that $ \m_{t+1}(y)=\sum_x\m_{t}(x)P(x,y)$.

Note that $\{\m^\a_t\}_{t\in \mathbb N}$ is a particular evolving measure with $\m_0=\a$ and also
$\{\m^\a_{t+t_0}\}_{t\in \mathbb N}$ is an  evolving measure, for fixed $t_0\in \mathbb N$.
 \item
 {\bf Squeezing measure and first hitting distribution:}
 as introduced in (\ref{mumu*}) a special role will be played by
the {\it squeezing-quasi-stationary measure} on $\cX$:
\be{defsq}
  \mu^{\mu^*}_t(y)=\sum_{z\in A}\m^*(z)P^t(z,y)=
  \begin{cases}
  \lambda^{t}\mu^{*}(y) & \text{ if }\; y \in A\\
  (1-\lambda^{t})\omega(y) & \text{ if }\; y\in G
  \end{cases}
\ee
where  the probability measure $\omega$ on $G$ is the {\it first hitting distribution} defined, for $y\in G$, by:
\be{defomega}
\omega(y)=\P\Big(X^{\m^*}_1=y\Big| X^{\m^*}_1\in G\Big)\equiv \frac{\sum_{z\in A}\m^*(z)P(z,y)}{1-\l}
\ee

Clearly  $\{\mu^{\mu^*}_t\}_{t\in \mathbb N}$ is an evolving measure.
\item
{\bf Separation ``distance":}
given two measures $\nu_1$ and $\nu_2$ on $\cX$  their
\emph{separation} is defined by
$$
sep(\nu_1,\nu_2):= \max_{y\in\cX}\big[1-\frac{\nu_1(y)}{\nu_2(y)} \big]
$$

 \end{itemize}
 
 \subsection{The local chain $\widetilde X_t$ on $A$}
\label{localch}
In this subsection we construct an ergodic Markov chain $\widetilde X_t$ on $A$, that we call 
the {\it local chain}. 

To describe the local behavior of the process $X_t$ on $A$, many different dynamics have been used in the literature.

The restriction of the transition matrix to the set $A$, $[P]_A$,  is a
sub-stochastic matrix, by adding to it a diagonal matrix $D$ with the escape probabilities,
$D(x,y)={\mathds 1}_{x=y}\sum_{z\in G}P(x,z)$,  one obtains the {\it reflected process} (see for instance
\cite{BG} and \cite{OV}) as a local dynamics. Another frequently used local process is the
{\it conditioned process}, defined by the original process $X_t$ on $\cX$
but conditioned to remain in $A$. This conditioned process has obviously a crucial role in the study of 
 the local behavior of the process $X_t$ before absorption in $G$. 
 However the main problem in dealing with it, for instance to estimate the hitting time to $G$, 
 is that this conditioned process  is no more a Markovian process.

We use here a different local chain  $\widetilde X_t$ constructed by means of the 
 right eigenvector  of $[P]_A$ corresponding to $\l$.
 This construction is related to the Doob h-transform of $[P]_A$ (see for instance \cite{LPW}).
This chain $\widetilde X_t$ is also related to the  ``reversed chain" in {Darroch-Seneta},introduced
in \cite{DS} while considering
the large time asymptotics. Our process $\widetilde X_t$ is the time reversal of this Darroch-Seneta
``reversed chain". 
\bigskip

The construction is the following:
 by the Perron-Frobenius theorem 
  there exists a unique non negative right eigenvector $\g$ of $[P]_A$ corresponding to $\l$, i.e.,
  \be{eqmqs1}
  [P]_A\g=\l \g\qquad \hbox{ with }\qquad (\m^*,\g)=1.
  \ee
This eigenvector is related to the asymptotic ratios
of the survival probabilities {(see eg \cite{CMS})}
$$
\lim_{t\to \infty}\frac{\P(\t^x_G>t)}{\P(\t^y_G>t)}=\frac{\g(x)}{\g(y)}\qquad x,y\in A.
$$

For any $x,y\in A$ define the stochastic matrix
\be{tildeP}
\widetilde P(x,y):=\frac{\g(y)}{\g(x)} \frac{P(x,y)}{\l}.
\ee
Let $\nu$ be its invariant measure
$$
\sum_{x\in A}\nu(x)\widetilde P(x,y)=\nu(y)=\sum_{x\in A}\nu(x) \frac{\g(y)}{\g(x)} \frac{P(x,y)}{\l}
$$
so that 

$$\g(x)=\frac{\nu(x)}{\m^*(x)}$$

For the chain $\widetilde X_t$ we define
$$\tilde s^x(t,y):=1-\frac{\widetilde P^t(x,y)}{\nu(y)}$$
$$\tilde s^x(t)= sep(\tilde \m^x_t,\nu)=\sup_{y\in A}\tilde s^x(t,y),\quad\tilde s(t):=\sup_{x\in A}\tilde s^x(t).$$
Note that $\tilde s^x(t)\in[0,1]$. Moreover, since $\widetilde P$ is a primitive matrix, it is well known (see for instance \cite{AD1}, Lemma 3.7) that
 $\tilde s(t)$ has 
 the sub-multiplicative property:
$$\tilde s(t+u)\le \tilde s(t)\tilde s(u).$$
This implies in particular an exponential decay in time of $\tilde s(t)$.

\bigskip

The relation between the local chain and the original chain $X_t$ on $\cX$
is given by the definition (\ref{tildeP}) and more generally by
\be{tildePt}
\widetilde P^t(x,y)=\frac{\g(y)}{\g(x)} \frac{P^t(x,y)}{\l^t}.
\ee
\bigskip

We can use  this relation 
to obtain a {\it rough estimate on the absorption time} $\t_G$.
We give here this simple calculation in order to point out the dependence on the initial distribution $\a$
of the distribution of $\t^{\a}_G$ by means 
 of  a  {\it time shift}.
 
 As it will be clear in what follows, it is natural to associate to every initial measure $\a$
the following measure $\tilde\a$ for the local chain $\widetilde X_t$:
$$
\tilde \a(x)=\frac{\a(x)\g(x)}{\sum_{y\in A}\a(y)\g(y)}.
$$

Indeed
 $$
\P(\t^{\a}_G>t)=\sum_{y\in A}\sum_{x\in A}\a(x) {P^t(x,y)}=$$
$$
\sum_{y\in A}\sum_{x\in A}\a(x) {\g(x)\l^t\m^*(y)\frac{\widetilde P^t(x,y)}{\nu(y)}}=
$$
$$
\sum_{y\in A}\sum_{x\in A}\a(x) \g(x)\l^t\m^*(y)(1-\tilde s^x(t,y)).$$
Since 
$\tilde s^x(t,y)\le \tilde s(t)$ we get
$$
\P(\t^{\a}_G>t)\ge (1-\tilde s(t)){\sum_{y\in A}}\sum_{x\in A}\a(x) \g(x)\l^t{\m^*(y)}=
\l^{t+\d_{\a}}(1-\tilde s(t))
$$
with
$$
\d_{\a}:=\log_\l\big(\sum_{x\in A}\a(x)\g(x)\big) 
$$

On the other side
we can consider the minimal  strong stationary time $\tilde \t^x_\nu$ such that
$$
\P(\widetilde X^x_t=y,\; \tilde \t^x_\nu=t)=\nu(y) \P(\tilde \t^x_\nu=t)
$$
 with
$$
\P(\tilde\t^x_\nu>t)=\tilde s^x(t).
$$
We have immediately
$$
\P(\t^{\a}_G>t)=\sum_{y\in A}\sum_{x\in A}\a(x) \g(x)\l^t\frac{\m^*(y)}{\nu(y)}\P(\widetilde X^x_t=y, \tilde\t^x\le t)+
$$
$$
\sum_{y\in A}\sum_{x\in A}\a(x) \g(x)\l^t\frac{\m^*(y)}{\nu(y)}\P(\widetilde X^x_t=y, \tilde\t^x> t)\le
$$
$$
\sum_{y\in A}\sum_{x\in A}\a(x) \g(x)\l^t\frac{\m^*(y)}{\nu(y)}\nu(y)\P(\tilde\t^x\le t)+
\sum_{y\in A}\sum_{x\in A}\a(x) \g(x)\l^t\frac{\m^*(y)}{\nu(y)}\P(\tilde\t^x>t)
$$
$$
=\l^{t+\d_{\a}}\Big[1+\tilde s^{\tilde\a}(t)\big(\sum_{y\in A} \frac{\m^*(y)}{\nu(y)}-1\big) \Big]
$$
with
$$
\tilde s^{\tilde\a}(t):=\frac{\sum_{x\in A}\a(x) \g(x) \tilde s^x(t)}{\sum_{x\in A}\a(x) \g(x)}
$$
Note that $\sum_{y\in A} \frac{\m^*(y)}{\nu(y)}\ge 1$. This quantity could be much larger that $1$
and so this estimate from above on the distribution of $\t^\a_G$ is quite rough due to the
factor $\big(\sum_{y\in A} \frac{\m^*(y)}{\nu(y)}-1\big)$. However, we have to note
that this factor is independent of time so that, for large $t$,  due to the exponential decay of $\tilde s(t)$,
and so
of $\tilde s^{\tilde\a}(t)$,  the estimate is not trivial. Similar results can be found in the literature.

\bigskip

 There are some interesting points to note  after this estimates, that will be important
 in our discussion especially for applications to metastability.
 \bi
\item
First of all we are able to consider arbitrary initial measures and to determine the effect of  the initial condition
on the distribution of the first hitting time to $G$.
Indeed we can associate to every initial measure $\a$ a corresponding time shift $\d_\a$. 
\item
We are interested in the application of  first hitting results to metastability.
In metastable situations the chain $\widetilde X_t$ has typically  a relaxation time much smaller than the mean absorption time
of the chain $X_t$.  This fast convergence
to equilibrium will be given by a fast exponential decay of the separation distance 
$$\tilde s^{\tilde \a}(t)=\max_{y\in A} \tilde s^{\tilde \a}(t,y).$$
Our control on the process $X_t$ with the local chain $\widetilde X_t$ given
by (\ref{tildePt}) is really strong  when looking at convergence to equilibrium in
separation distance, see Proposition \ref{smu0}.
This implies that we can  use the good convergence to equilibrium of $\widetilde X_t$ to obtain better estimates 
on the absorption time
of the chain $X_t$.
\ei

 \subsection{Main results}
\label{mainres} 
 
 We first extend the notion of Strong Stationary Time (SST)  to {\it strong time w.r.t. evolving measures}, different from
the stationary one, with the following.
 \bd{smureft}
For any initial distribution $\a$   and for any  evolving measure $\{\m_t\}_{t\in \mathbb N}$ on $\cX$, 
we call    \emph {strong  time w.r.t. $\{\m_t\}_{t\in \mathbb N}$} a randomized stopping time, $\t^{\a}_\m$, for $X^{\a}_t$ such that for any $y\in \cX$ we have
 \be{defSTwrt}
 \P(X^{\a}_t=y,\, \t^{\a}_\m=t)=\m_t(y)\P(\t^{\a}_\m=t)
 \ee
\ed
Note that:
\be{STwrt1}
 \P(X^{\a}_t=y,\, \t^{\a}_\m\le t)=\m_t(y)\P(\t^{\a}_\m\le t).
\ee
Indeed 
$$
\P(X^{\a}_t=y,\, \t^{\a}_\m\le t)=\sum_{u=0}^t\P(X^{\a}_t=y,\, \t^{\a}_\m=u)=
\sum_{u=0}^t\sum_{z\in\cX}\P(X^{\a}_u=z,\, \t^{\a}_\m=u)\P(X^{z}_{t-u}=y)=$$
$$
\sum_{u=0}^t\sum_{z\in\cX}\P( \t^{\a}_\m=u)\m_u(z)P^{t-u}(z,y)=\m_t(y)\P( \t^{\a}_\m\le t)
$$
For these strong times we have a result similar to what is proved in \cite{AF} for strong stationary
times and their relation with separation distance between the evolution and the stationary measure.
\bt{main0}
  For any initial distribution $\a$  and
  for any reference evolving measure $\{\m_t\}_{t\in \mathbb N}$ and any strong time  $\t^{\a}_\m$ w.r.t. $\{\m_t\}_{t\in \mathbb N}$   we have
 $$
 \P(\t^{\a}_\m>t)\ge sep(\m^{\a}_t, \m_t).
 $$ 
 Moreover there exists a \emph{minimal} strong time w.r.t. $\{\m_t\}_{t\in \mathbb N}$  such that
 $$
 \P(\t^{\a}_\m>t)= sep(\m^{\a}_t, \m_t).
 $$
\et

The proof of Theorem \ref{main0} is given in Section \ref{pf_main0}. 
 \bigskip
 
By using the time shift associated to the initial measure $\a$, obtained in
the rough estimate on the absorption time given in the previous section,
we are going to define the  {\it sqeezed-quasi-stationary reference measure} $\r_t$.

Recall that for any initial distribution $\a$ we define the \emph{time shift} 
$$
\d_{\a}:=\log_\l\big(\sum_{x\in A}\a(x)\g(x)\big) 
$$
with $\g$ defined in (\ref{eqmqs1}).
\bd{timesh}
The following {\emph {reference evolving measure }} $\{\rho_t\}_{t\in \mathbb N}$
depending on $\a$:
$$
\r_t:= \m^{\m^*}_{t+\d_{\a}}
$$
is a probability measure if $t+\d_{\a}\ge 0$.
\ed
Moreover for any $y\in A$ and  $t+\d_{\a}\ge 0$ define: 
$$
s^{\a}(t,y)=1-\frac{\m^{\a}_t(y)}{\r_t(y)}=1-\frac{\m^{\a}_t(y)}{\m^{\m^*}_{t+\d_{\a}}(y)}
$$
and 
$$
\tilde s^{\tilde\a}(t,y)=1-\frac{\sum_{x\in A}\tilde\a(x) \widetilde P^t(x,y)}{\nu(y)}
$$

\bp{smu0}
For any $y\in A$  we have: 
$$
s^{\a}(t,y)=\tilde s^{\tilde\a}(t,y)
$$
\ep
The proof is immediate since
$$
s^{\a}(t,y)=1-\frac{\sum_{x\in A}\a(x) P^t(x,y)}{\l^{t+\d_{\a}}\m^*(y)}=
1-\frac{\sum_{x\in A}\a(x)\g(x)\l^t\frac{\m^*(y)}{\nu(y)} \widetilde P^t(x,y)}{\l^{t+\d_{\a}}\m^*(y)}=
$$
$$
1-\frac{\sum_{x\in A}\a(x)\g(x)(1-\tilde s^x(t,y))}{\sum_{x\in A}{\a(x)\g(x)}}=\tilde s^{\tilde\a}(t,y).
$$
\eproof

\bigskip

By Theorem \ref{main0} from the separation $s^{\a}(t)$ we can define a minimal strong
time w.r.t. the reference measure $\r_t$, say $\t^\a_\r$, such that
$$
\P(X^\a_t=y, \, \t^\a_\r=t)=\r_t(y)\P(\t^\a_\r=t)\qquad \hbox{ with }\qquad \P(\t^\a_\r>t)=s^{\a}(t).
$$
With a simple argument we have
$$
\P(\t^\a_G>t)=\sum_{s\le t}\sum_{y\in A}\P(X^\a_t=y, \, \t^\a_\r=s)+\P(\t^\a_G>t,\;\t^\a_\r>t)=
$$
\be{stro}
\l^{t+\d_\a}\big(1-s^\a(t)\big)+\P(\t^\a_G>t,\;\t^\a_\r>t).
\ee
If $s^\a(t)$ decays in time faster than $\l^{t+\d_\a}$, we can obtain from (\ref{stro})
good estimates from above and from below, i.e., we get
$$
\P(\t^\a_G>t)=\l^{t+\d_\a}\big(1+o(1)\big).
$$
Notice that, while $\tilde s^{\tilde \alpha}(t)$  decays exponentially in time, in general we 
don't have such a good long time behavior for
$$
s^{\a}(t)=\max_{y\in\cX}s^{\a}(t,y)=\max_{y\in A}\tilde s^{\tilde\a}(t,y)\vee \max_{y\in G}s^{\a}(t,y).
$$
Indeed, even in metastable situations, where we expect a decay of $\tilde s^{\tilde\a}(t)$
faster do than $\l^{t+\d_\a}$, in general we cannot control the term $\max_{y\in G}s^{\a}(t,y)$.

To solve this problem, we
 define a new random time by looking at the conditioned process by means of the separation
$\tilde s^{\tilde\a}(t)$ instead of $s^{\a}(t)$.

\bd{csqst}
 For any initial distribution $\a$ on $A$ we call    {\emph {conditionally strong quasi stationary
 time (CSQST)}} a randomized stopping time $\t^{\a}_*$ for $X^{\a}_t$ such that for any $y\in A$ we have
  $$
  \P(X^{\a}_t=y,\, \t^{\a}_*=t\big| t<\t^{\a}_G)=\m^*(y)\P(\t^{\a}_*=t\big| t<\t^{\a}_G)
  $$
 which is equivalent to 
 \be{defCSQST}
 \P(X^{\a}_t=y,\, \t^{\a}_*=t)=\m^*(y)\P(\t^{\a}_*=t<\t^{\a}_G)
 \ee
  \ed
We note that the analogous of the equation (\ref{STwrt1}) holding for strong times,
does not hold for CSQST. Due to the conditioning, we have:
\be{defCSQSTle}
 \P(X^{\a}_t=y,\, \t^{\a}_*\le t)= \m^*(y)\sum_{u\le t}\l^{t-u}\P(\t^{\a}_*=u<\t^{\a}_G)\not=
  \m^*(y)\P(\t^{\a}_*\le t <\t^{\a}_G)
 \ee
 indeed
 $$
\sum_{u\le t}\sum_{z\in A}\P(X^{\a}_u=z,\, \t^{\a}_*=u)P^{t-u}(z,y) =\sum_{u\le t}\sum_{z\in A}\m^*(z)\P(\t^{\a}_*=u<\t^{\a}_G)P^{t-u}(z,y)=
$$
$$
\m^*(y)\sum_{u\le t}\l^{t-u}\P(\t^{\a}_*=u<\t^{\a}_G).$$
This remark actually suggests a new notion of minimality as given in the following.

\bt{main1}
  For any initial distribution $\a$ on $A$ and
 for any  $\t^{\a}_*$   { {conditionally strong quasi stationary
 time (CSQST)}} for $X^{\a}_t$ and for all $t\ge 0$ we have
 $$
 \sum_{u\le t}\l^{-u}\P(\t^{\a}_*=u<\t^{\a}_G)\le \l^{\d_{\a}}(1-\tilde s^{\tilde\a}(t)).
 $$
 Moreover there exists a \emph{minimal conditionally strong quasi stationary
 time} $\t^{\a}_*$  such that
 $$
 \sum_{u\le t}\l^{-u}\P(\t^{\a}_*=u<\t^{\a}_G)= \l^{\d_{\a}}(1-\tilde s^{\tilde\a}(t)).
 $$
with 
$$\P(\t^{\a}_*=t<\t^{\a}_G)= \l^{t+\d_{\a}}(\tilde s^{\tilde\a}(t-1)-\tilde s^{\tilde\a}(t)).$$
 \et

 Note that in particular for a minimal conditionally strong quasi stationary
 time we have
$$\P(\t^{\a}_*>t, \; \t^{\a}_*<\t^{\a}_G)= \sum_{u>t}\l^{u+\d_{\a}}(\tilde s^{\tilde\a}(u-1)-\tilde s^{\tilde\a}(u))
\le \l^{t+\d_\a}\tilde s^{\tilde\a}(t).$$

The interest of this  {minimal} conditionally strong quasi stationary
 time is given by the following:
 \bt{main2}
 For any initial distribution $\a$ on $A$,  if  $\t^{\a}_*$   is a minimal conditionally strong quasi stationary
 time and $t+\d_{\a}\ge 0$ we have
 $$
 \P\Big(\t^{\a}_G>t\Big)=\l^{t+\d_{\a}}(1-\tilde s^{\tilde\a}(t))+\P\Big( \t^{\a}_{*,G}> t\Big)
 $$
 with $\t^{\a}_{*,G}={\t^{\a}_G\wedge \t^{\a}_*}$.
 
 Moreover for any $y\in G$ we have
 \be{uscita}
 \P\Big(X^{\a}_{\t^{\a}_G}=y\Big)=\P\Big( \t^{\a}_G< \t^{\a}_*, \; X^{\a}_{\t^{\a}_G}=y\Big)+
 \omega(y)\P\Big( \t^{\a}_G> \t^{\a}_*\Big).
\ee
 \et 
  \bigskip
  
  This theorem provides a quantitative control on the convergence  to an exponential distribution
  for the hitting time $\t_G$ and  on the exit distribution.
  Note that with this CSQST, $\t^{\a}_{*}$, we are obtaining  conditioning benefits 
  without  explicitly using the conditioned process.

  As far as the exit distribution is concerned 
   in the metastable case  the quantity
 $$\P\Big( \t^{\a}_G< \t^{\a}_*\Big)
 =:\e$$
 should be  small. This is the case in which
 the results of Theorem \ref{main2} are relevant.
I, indeed, $\e$ small  implies that
  the distribution of the
 first hitting to $G$ is well approximated by the measure $\omega$ since
 equation (\ref{uscita}) gives
 $$
 \omega(y)(1-\e)\ge \P\Big(X^{\a}_{\t^{\a}_G}=y\Big)\le \e+\omega(y)(1-\e)
 $$

 
 \section{Strong time w.r.t.   evolving measures}
 \label{pf_main0}
 
 In this section we extend the results obtained in \cite{AF}, relating strong stationary times
 and separation distance, to strong times w.r.t. evolving measures.
 The ideas of the proof are simple.
 
  \bigskip
 
\subsection{Proof of Theorem \ref{main0}}
\label{s1}

For any $y\in\cX$ and $t\ge 0$ and any evolving measure $\m_t$, define
$$s^{\a}_\m(t,y) := 1-\frac{\m^{\a}_t(y)}{\m_t(y)}, \quad sep (\m^{\a}_t,\m_t):= s^{\a}_\m(t)=\sup_{y\in\cX}s^{\a}_\m(t,y).$$
Let $\t^{\a}$ be a strong time w.r.t.$\m_t$ then
$$
\m^{\a}_t(y)\ge\P(\t^{\a}\le t, X^{\a}_t=y)
$$
so that 
for any $t\ge 0$ and $y\in\cX$ 
$$
1-s^{\a}_\m(t,y)=\frac{\m^{\a}_t(y)}{\m_t(y)}\ge\frac{\P(\t^{\a}\le t, X^{\a}_t=y)}{\m_t(y)}=\P(\t^{\a}\le t)
$$
and so
$$\P(\t^{\a}>t)\ge s^{\a}_\m(t).$$
On the other side starting from the separation $s^{\a}_\m(t)=sep(\m^{\a}_t, \m_t)$ we can define a
minimal strong  time w.r.t. $\m_t$ as follows:
note that $s^{\a}_\m(t) \in [0,1]$ and it is a decreasing function of $t$, say $s^{\a}_\m(t+1)\le s^{\a}_\m(t)$. Define $s^{\a}_\m(-1):=1$ and
$$
\s_t(y):=\m_t(y) [s^{\a}(t-1)-s^{\a}(t)]\qquad \theta_t(y):=\m_t(y) [s^{\a}_\m(t-1)-s^{\a}_\m(t,y)]
$$
We have for any $y\in\cX$ and $t\ge 0$
$$0\le\s_t(y)\le\theta_t(y)$$
and more precisely
$$
\theta_t(y)-\s_t(y)=\m_t(y)[s^{\a}_\m(t)-s^{\a}_\m(t,y)]
$$
so that the vectors $\s_t$ and $\theta_t$ satisfy the {iterative equation}
\be{ric_st}
(\theta_t-\s_t)P=\theta_{t+1} \qquad \forall t\ge 0 
\ee
 with $ \theta_0=\a,$ and $ \s_t= \Big(\min_{z\in \cX}\frac{\theta_t(z)}{\m_t(z)}\Big)\m_t$ for all $ t\ge 0 $.

Define a randomized stopping time $\t^{\a}$ by imposing
\be{def_tau}
\P\Big(\t^{\a}=t\Big| \t^{\a}\ge t, X^{\a}_t=y, X_s, s<t\Big)=
\frac{\s_t(y)}{\theta_t(y)}=\frac{s^{\a}_\m(t-1)-s^{\a}_\m(t)}{s^{\a}_\m(t-1)-s^{\a}_\m(t,y)}
\ee
It is easy to prove by induction that for any $t\ge 0$
\be{ts1}
\P(\t^{\a}=t, X^{\a}_t=y)=\s_t(y), \quad \P(\t^{\a}\ge t, X^{\a}_t=y)=\theta_t(y)
\ee
since also these probabilities satisfy the iterative equation (\ref{ric_st}).
Indeed if (\ref{ts1}) holds for $t$ then by (\ref{ric_st}) we obtain the statement for 
$\theta_{t+1}$ and by (\ref{def_tau}) the same for $\s_{t+1}$.
We can immediately
conclude that $\t^{\a}$ is a  strong  time w.r.t. $\m_t$ with $\P(\t^{\a}>t)=s^{\a}_\m(t)$.
Thus, it is minimal.

\bigskip
\subsection{Construction of the strong time w.r.t. $\{\m_t\}_{t\in \mathbb N}$  with an auxiliary chain}
\label{s2}
We give here a   construction of the minimal
strong time  $\t^{\a}$ inspired by \cite{AD2,DM09}.  We define an auxiliary chain
so that the strong time can be seen as an hitting time for this new process. 

Consider the initial distribution ${\a}$ as a parameter and define an auxiliary process $Y^{\a}_{t}$ with state
space $\mathcal{Y}:=\mathcal{X}\times\{0,1\}$, so that on $\{0\}$ the process is like
$X^{\a}_t$ but with a rate jump to $\{1\}$ given by $J^{\a}$.

More precisely for every $z\in\cX$ define the function 
\begin{equation}
J^{\a}(t,z):=\frac{s^{\a}_\m(t-1)-s^{\a}_\m(t)}{s^{\a}_\m(t-1)-s^{\a}_\m(t,z)},\label{eq:defJ}
\end{equation}
with the convention $0/0=0$.
B, by the monotonicity of  $s^{\a}(t)$, we have $J^{\a}(t,z)\in[0,1]$
for any $ z\in\cX$ and any $t$ . 

Consider the following time dependent transition probabilities
for the process $Y_t^{\a}$:
$$Q^{\a}_{(y,0),(z,0)}=P(y,z)\Big(1-J^{\a}(t,z)\Big), \quad 
Q^{\a}_{(y,0),(z,1)}=P(y,z)J^{\a}(t,z),\quad
Q^{\a}_{(y,1),(z,e)}=P(y,z)\d_{1,e}.
$$
Note that the marginal distribution of $Y_t^{\a}$ on $\cX$ corresponds to the distribution
of $X^{\a}_t$ so that we can study each event defined for the process $X^{\a}_t$ in terms
of set of paths of the process $Y^{\a}_t$. For this reason, with an  abuse of notation, we denote with the same symbol
$\P$ the probability of events defined in terms of  the process  $Y^{\a}_t$.
Consider the hitting time:
\[
\tau^{\a}_\mathtt{1}:=
\tau^{\a}_{\cX \times \{1\}}=\min\{t\ge0\:;\: Y_{t}^{\a}=(y,1)\:\text{for some }y\in\mathcal{X}\},
\]
We want to show that $\tau^{\a}_\mathtt{1}$ is a minimal strong time w.r.t. the evolving measure $\{\m_t\}_{t\in \mathbb N}$, i.e., 
\be{pa1}
 \P(X^{\a}_t=y,\, \t^{\a}_1=t)=\m_t(y)\P(\t^{\a}_1=t)=\m_t (y)\Big( s^{\a}_\m(t-1)-s^{\a}_\m(t)\Big)
\ee

We proceed by induction on $t$.
For $t=0$, by definition of  $Y^{\a}_{t}$, we have 
$$\P(X^{\a}_0=y,\, \t^{\a}_1=0)=
\P(Y^{\a}_0=(y,1))=\a(y)J^{\a}(0,y)=$$
$$\a(y)\frac{1-s^{\a}_\m(0)}{1-s^{\a}_\m(t,y)}=\m_0(y)(1-s^{\a}_\m(0)).$$

To prove the induction step we use the following:
\bl{l001}
If for any $u\le t$ we have
$$
\P\big( X^{\a}_u=y\;|\t^{\a}_1=u\big)={\m}_u(y)
$$
then
$$
\P\big( Y^{\a}_t=(z,1)\big)={\m}_t(z)\P(\t^{\a}_1\le t)
$$
\el
\bpr
$$
\P\big( Y^{\a}_t=(z,1)\big)=\sum_{u\le t}\sum_{y\in\cX}\P\big( Y^{\a}_t=(z,1)|Y^{\a}_u=(y,1)\big)
\P\big(\t^{\a}_1=u, X^{\a}_u=y\big)=$$
$$\sum_{u\le t}\sum_{y\in\cX}P^{t-u}(y,z){\m}_u(y)
\P(\t^{\a}_1= u)=
\m_t(z)\P(\t^{\a}_1\le t)
$$
\epr

Suppose now that (\ref{pa1}) holds for $u\le t$.
By using then Lemma  \ref{l001} we get
$$
\P\big( X^{\a}_{t+1}=y, \t^{\a}_1={t+1} \big)=\sum_{z\in\cX}\P\big(Y^{\a}_{t+1}=(y,1)|Y^{\a}_{t}=(z,0)\big)\P\big(Y^{\a}_{t}=(z,0)\big)=
$$
$$
\sum_{z\in\cX}P(z,y)J^{\a}({t+1},y)\big[ \m^{\a}_{t}(z)-\P(Y^{\a}_{t}=(z,1)) \big]=$$
$$
J^{\a}({t+1},y)\big[\m^{\a}_{t+1}(y)-\sum_{z\in\cX}{\m}_{t}(z)P(z,y)\P(\t^{\a}_1\le t)\big]=
$$
$$
J^{\a}({t+1},y){\m}_{{t+1}}(y)\big[1-s^{\a}_\m({t+1},y)-(1-\P(\t^{\a}_1>t)\big]=$$
$$
{\m}_{{t+1}}(y)\frac{s^{\a}_\m(t)-s^{\a}_\m({t+1})}{s^{\a}_\m(t)-s^{\a}_\m({t+1},y)}\big[s^{\a}_\m(t)-s^{\a}_\m({t+1},y)\big]
$$
and summing on $y$ we get
$\P(\t^{\a}_1={t+1})=s^{\a}_\m(t)-s^{\a}_\m({t+1})=\P(\t^{\a}_1>t)-s^{\a}_\m({t+1})$ so that $\P(\t^{\a}_1>{t+1})=s^{\a}_\m({t+1})$ and
$$
\P\big( X^{\a}_{t+1}=y, \t^{\a}_1={t+1} \big)={\m}_{{t+1}}(y)\P(\t^{\a}_1={t+1}).\hskip3cm \qed
$$

\section{ Conditionally strong quasi stationary times (CSQST)}
\label{s_main1}
In this section we prove Theorem \ref{main1}.
The main idea is the following. As noted after Proposition \ref{smu0},
 given
a reference evolving measure $\r_t$, obtained by   the squeezing-quasi-stationary measure with time shift, we know how to construct
from the corresponding  separation distance a  minimal strong time w.r.t. $\r_t$.
Proposition \ref{smu0} opens the way to the construction of a faster strong time (with finite moments) if we take the supremum of $s^\alpha(t,y)$ only for $y$ in $A$, since this quantity coincides with $\tilde s^{\tilde \alpha}(t)$ that, being the separation from stationarity for the process $\widetilde X$, decays exponentially in time.
The idea is then to define a new random time $\t^\alpha_*$ by using $\tilde s^{\tilde\a}(t)$ instead of
$ s^{\a}(t)$, following the
construction given in the proof of Theorem \ref{main0}. This time is not strong w.r.t. $\r_t$, but
it works like  a strong time when looking at the process conditioned to $A$. This construction gives
 a conditionally strong-quasi-stationary time without working directly with the conditioned
process. 

We first prove that if $\t^{\a}_*$ is a CSQST, i.e., if satisfies
 $$
 \P(X^{\a}_t=y,\, \t^{\a}_*=t)=\m^*(y)\P(\t^{\a}_*=t<\t^{\a}_G)
$$
 then 
 for all $t\ge 0$ we have
 \be{bo1}
 \sum_{u\le t}\l^{-u}\P(\t^{\a}_*=u<\t^{\a}_G)\le \l^{\d_{\a}}(1-\tilde s^{\tilde\a}(t)).
 \ee
Indeed for any $y\in A$ we have
\begin{eqnarray*}
\m^{\a}_t(y)
&\ge& 
\P(\t^{\a}_*\le t, X^{\a}_t=y)=
 \sum_{u\le t}\sum_{z\in A}\P(\t^{\a}_*=u, X^{\a}_u=z)P^{t-u}(z,y)=\\
&&\l^t
 \sum_{u\le t}\l^{-u}\P(\t^{\a}_*=u<\t^{\a}_G)\m^*(y)
\end{eqnarray*}
so that 
$$
\frac{\m^{\a}_t(y)}{\l^t\m^*(y)}=\l^{\d_{\a}}(1-s^{\a}(t,y))\ge  \sum_{u\le t}\l^{-u}\P(\t^{\a}_*=u<\t^{\a}_G)
$$
since this holds for any $y\in A$ and  we have $s^{\a}(t,y)=\tilde s^{\tilde\a}(t,y)$ for every $y\in A$  then (\ref{bo1}) holds.

We define now a random time
$\t^{\a}_*$ which is not strong w.r.t. the reference evolving measure  
$\r_t\equiv \m^{\m^*}_{t+\d_{\a}}$
on the
hole space $\cX$ but which is constructed with similar ideas, by using
Proposition \ref{smu0}, by means of
the separation $\tilde s^{\tilde \a}$ in the following way.
Define
$$
\s_t(y)={\mathds 1}_{y\in A}\;\r_t(y)\big(\tilde s^{\tilde\a}(t-1) -\tilde s^{\tilde\a}(t) \big)
$$
$$
\theta_t(y)={\mathds 1}_{y\in A}\;\r_t(y)\big(\tilde s^{\tilde\a}(t-1) -\tilde s^{\tilde\a}(t,y) \big)
$$
Then for any $y\in A$ we still can define $\t^{\a}_*$ such that
$$
\P(X^{\a}_t=y, \t^{\a}_*=t)=\s_t(y)=\r_t(y)\big(\tilde s^{\tilde\a}(t-1) -\tilde s^{\tilde\a}(t) \big)=
$$
$$
\m^*(y)\P(\t^{\a}_*=t<\t^{\a}_G)
$$
with
$$
\P(\t^{\a}_*=t<\t^{\a}_G)=\l^{t+\d_{\a}}\big(\tilde s^{\tilde\a}(t-1) -\tilde s^{\tilde\a}(t) \big)
$$
$$
\P(\t^{\a}_*=t\ge\t^{\a}_G)=0
$$

and
$$
\P(\t^{\a}_*=+\infty)=1-\sum_{t\ge 0}\l^{t+\d_{\a}}\big(\tilde s^{\tilde\a}(t-1) -\tilde s^{\tilde\a}(t) \big)=
\P(\t^{\a}_G<\t^{\a}_*)>0
$$
For such a $\t^{\a}_*$ we have:
\bp{ptstar}
 $\t^{\a}_*$ is a conditionally strong quasi-stationary  time, i.e.,
$$
\P\Big(X^{\a}_t=y, \t^{\a}_*=t\Big|\t^{\a}_G>t\Big)=
\m^*(y) \P\Big( \t^{\a}_*=t\Big|\t^{\a}_G>t\Big)
$$
\ep
\bigskip

Indeed
$$
\P\Big(X^{\a}_t=y, \t^{\a}_*=t\Big|\t^{\a}_G>t\Big)=
\frac{\P\Big(X^{\a}_t=y, \t^{\a}_*=t<\t^{\a}_G\Big)}
{\P\Big(\t^{\a}_G>t\Big)}=
$$
$$
\frac{\m^*(y)\l^{t+\d_{\a}}(\tilde s_{\a}(t-1)-\tilde s_{\a}(t))}{\P\Big(\t^{\a}_G>t\Big)}=
\m^*(y)\frac{\P\Big(\t^\a_*=t<\t^{\a}_G\Big)}{\P\Big(\t^{\a}_G>t\Big)}\qquad \qed
$$

\section{ Representation formula for $\t^{\a}_G$ with $\t^\a_*$ }
\label{s_main2}
In this section we prove Theorem \ref{main2}.
We first prove that
$$
\P(\t^{\a}_G>t)=\l^{t+\d_{\a}}(1-\tilde s^{\tilde\a}(t))+\P( \t^{\a}_{*,G}> t)
$$

Indeed we have
$$
\P(\t^{\a}_G>t)=\P(\t^{\a}_G>t, \t^{\a}_*\le t)+\P(\t^{\a}_G>t, \t^{\a}_*> t)=$$
$$
\sum_{y\in A} \P(X^{\a}_t=y, \t^{\a}_*\le t)+\P({\t^{\a}_G\wedge \t^{\a}_*}> t)=$$
$$
\sum_{y\in A}\sum_{z\in A}\sum_{u=0}^t\P(X^{\a}_u=z, \t^{\a}_*=u, X^{\a}_t=y)+\P( {\t^{\a}_{*,G}}> t)=
$$
$$
\sum_{y\in A}\sum_{z\in A}\sum_{u=0}^t\m^*(z)\l^{u+\d_{\a}}(\tilde s^{\tilde\a}(u-1)-\tilde s^{\tilde\a}(u))P^{t-u}(z,y)+\P( \t^{\a}_{*,G}> t)=
$$
$$
\sum_{y\in A}\sum_{u=0}^t\m^*(y)\l^{u+\d_{\a}+t-u}(\tilde s^{\tilde\a}(u-1)-\tilde s^{\tilde\a}(u))+\P( \t^{\a}_{*,G}> t)=$$

$$\l^{t+\d_{\a}}(1-\tilde s^{\tilde\a}(t))+\P( \t^{\a}_{*,G}> t).
$$
Moreover for any $y\in G$ we have
 $$
 \P\Big(X^{\a}_{\t^{\a}_G}=y\Big)=\P\Big( \t^{\a}_G< \t^{\a}_*, \; X^{\a}_{\t^{\a}_G}=y\Big)+
\P\Big( \t^{\a}_G> \t^{\a}_*, \; X^{\a}_{\t^{\a}_G}=y\Big) 
$$
The second term in the r.h.s. can be written as
$$
\sum_{t=0}^\infty\sum_{z\in A}\P\Big( \t^{\a}_G> t=\t^{\a}_*, \: X^{\a}_t=z\Big)\P\Big(X^z_{\t^{z}_G}=y\Big)
=\sum_{t=0}^\infty\sum_{z\in A}\m^*(z)\P\Big( \t^{\a}_G> t=\t^{\a}_*\Big)\P\Big(X^z_{\t^{z}_G}=y\Big)$$

$$=\P\Big( \t^{\a}_G>\t^{\a}_*\Big)\sum_{z\in A}\m^*(z)\sum_{u=0}^{\infty}\sum_{w\in A}\P(X^z_u=w, \; \t^z_G=u+1)=
\omega(y)\P\Big( \t^{\a}_G> \t^{\a}_*\Big)
$$
so that (\ref{uscita}) holds.

\section{Concluding remarks and future perspectives}
\label{s5}

In this paper we describe the relation between rarity and exponentiality with the help of a new class
of strong times.
We give an exact representation formula (Theorem\ref{main2}) that provides probabilistic interpretations
for the leading exponential term as well as for the error term.

Our setting is completely general: we do not need reversibility  and we only assume that $[P]_A$ is
a primitive matrix. Our representation formula applies to any initial state $\a$.
To our knowledge, no other result is so general and so transparent about the role of the starting state.
As discussed in the introduction, in the literature many results about exponentiality of the first hitting
time are obtained with renewal arguments based on the idea of recurrence to a ``basin of attraction" of
the metastable state. 
The control that we have in this paper on the role of the starting state $\a$ is such that we can obtain estimates 
on the distribution of $\t^\a_G$ without recurrence on a particular set but converging to a particular evolving measure,
that depends on the initial distribution $\a$, without error propagation. 
Indeed, we associate to each starting distribution $\a$ a time shift $\d_\a$, in such a way that the dependence on the initial 
distribution of the distribution of $\t^\a_G$ is described in terms of this time shift.
The evolving measure associated to the starting distribution $\a$ with this time shift is a probability
measure for every $t\ge 0\vee ( - \delta_\alpha)$
The set of states with $\delta_\alpha <0$ can be seen as metastable basin.

The main novelty of this paper is the introduction of a new language to describe the hitting of a set in terms of strong times.
Under very general conditions, the distribution of Conditionally strong quasi stationary times (see Def. \ref{csqst}) has a good asymptotic behavior.
In many physical applications however, one is more interested in the short time behavior of the process, and our notion of strong time w.r.t. other evolving measures (see Def \ref{smureft}) may give interesting estimates for such small times. 


Most of the  bounds of the error term in the exponential approximation of hitting times known in the literature  are function of the ratio
between a ``mean local relaxation time" and the mean hitting time. At heuristic level this time-scale
comparison is a very popular characterization of metastability. One of our strongest motivations
has been to give a rigorous base to this idea and to give a general characterization of metastability in terms of a time comparison.
This is still the first point in our agenda. 
In Theorem \ref{main2} these two different time scales are given by the times $\t^\a_{*,G}$ (that plays the
role of local relaxation time) and $\t^\a_G$.

The usability of our representation formula in Theorem \ref{main2} to get explicit error bounds relies on the
possibility to estimate $\tilde s^{\tilde\a}(t)$. At first glance, this task seems rather difficult because this
quantity is defined in terms of the matrix $\widetilde P$ and of the eigenvectors $\m^*$ and $\g$ of $[P]_A$.
Moreover, generally speaking, separation is not the most manageable notion of distance between measures.
However, since $\tilde s^{\tilde\a}(t)\equiv sep(\tilde\m^{\tilde\a}_t,\nu)$ is the separation for the chain
$\widetilde X_t$, it is positive; most important, it is bounded above by $\tilde s(t):=\sup_{\tilde\a}\tilde s^{\tilde\a}(t)$
which is submultiplicative.
Therefore, it is sufficient to find a time $R$ for which $\tilde s(t)$ is bounded above by a constant $c$ smaller than $1$
to get an exponential bound like $c^{t/R}$.
Useful inequalities that relate separation from stationarity are known (see e.g. \cite{AD1}) and can be used to find such a bound.

In order to control the effect of the initial distribution $\a$ on the distribution of $\t^\a_G$, a crucial tool turns out
to be the local chain $\widetilde X_t$. 
The main feature of this local chain is given by Proposition \ref{smu0}, which allows to see that 
 $\sup_{y\in A}s^\a(t,y)=\tilde s^{\tilde\a}(t)$ decays exponentially uniformly
in $\a$. 
Proposition \ref{smu0} also allows to compute $\tilde s^{\tilde\a}(t)$ without computing $\widetilde P, \m^*$ and $\g$.
In metastable situations we expect that this exponential decay is much faster that the decay of
$s^\a(t)=\sup_{y\in \cX}s^\a(t,y)$.
In the strong-time language, this means that $\t^\a_\r$ is slower than $\t^\a_*$ for it triggers the 
arrival to $\r_t$ also for the points in $G$.
By using the notion of CSQST with the representation formula of Theorem \ref{main2} we can use
the fast decay of $\tilde s(t)$ in order to control the distribution of $\t^\a_G$ by means of the distribution
of $\t^\a_{*,G}$.
The exponential behavior emerges when 
the term $\P\Big( \t^{\a}_{*,G}> t\Big)$
has a decay strictly faster than $\l^{t+\d_\a}$,
a sort of time comparison that may be used 
to characterize metastability.
  
  The statement of Theorem \ref{main2}
  has a strong analogy with the description
  of
  metastability in terms  of  recurrence \cite{FMNS}, \cite{FMNSS}.
  In the simple case of recurrence to a single state $x_0$, the main metastability hypothesis was
  on the decay in time of the quantity $sup_{x\in\cX}\P(\t^x_{x_0\cup G}>t)$, here replaced by
  a decay of $\P\Big( \t^{\a}_{*,G}> t\Big)$.
  Moreover,  an analogous of the auxiliary chain given in subsection
  \ref{s2} can be defined to see  $\t^{\a}_*$  as a hitting time. In this way we expect that exponential
  estimates from above can be obtained for the conditioned probability $\P\Big( \t^{\a}_{*}> t\Big|\t^{\a}_G> t \Big)$,
   with arguments similar to those used to estimate $sup_{x\in\cX}\P(\t^x_{x_0\cup G}>t)$ in some examples,
   see for instance \cite{FMNS}.
\bigskip

{\bf Acknowledgments: }

We thank Amine Asselah, Nils Berglund,  Pietro Caputo, Frank den Hollander, Roberto Fernandez and Alexandre Gaudilli\`ere for many fruitful  discussions.
 This work was partially supported by the A*MIDEX project (n. ANR-11-IDEX-0001-02) funded by the  ``Investissements d'Avenir" French Government program, managed by the French National Research Agency (ANR).



\begin{thebibliography}{99}
\bibitem{AG00}
 M.Abadi,  A.Galves ``Inequalities for the occurrence times of rare events in mixing processes. 
 The state of the art"
\emph{ Markov Process. Relat. Fields} {\bf 7}, 97--112 (2001).


\bibitem{A}
 D. Aldous,  ``Markov chains with almost exponential hitting times''
\emph{ Sto.Proc.Appl} {\bf 13}, 305--310 (1982).

\bibitem{AB1}
D.\ Aldous, M.\ Brown,
``Inequalities for rare events in time reversible Markov chains I'',
in \emph{Stochastic Inequalities}, M. Shaked and Y.L. Tong eds., pp. 1--16,
Lecture Notes of the Institute of Mathematical Statistics, vol. 22 (1992).

\bibitem{AB2}
 D.\ Aldous, M.\ Brown,
`` Inequalities for rare events in time reversible Markov chains II'',
\emph{ Sto.Proc.Appl} {\bf 44}, 15-25 (1993).

\bibitem{AD2}
 D.\ Aldous, P.Diaconis,
``Shuffling cards and stopping times'',
\emph{ Amer. Math. Monthly} {\bf 93}, 333-348 (1986).

\bibitem{AD1}
 D.\ Aldous, P.Diaconis,
``Strong uniform times and finite random walks I'',
\emph{ Adv. in Appl. Math.} {\bf 8}, 66-97 (1987).

\bibitem{AF}
 D.\ Aldous, J.A.Fill,
``Reversible Markov Chains and Random Walks on Graphs'',
\emph{Unfinished monograph, {2002}, recompiled 2014}, available 
      at \url{http://www.stat.berkeley.edu/~aldous/RWG/book.html} 


\bibitem{BG}
 A.\ Bianchi, A.\ Gaudilliere,
 ``Metastable states, quasi-stationary and soft masures, mixing time asymptotics via variational principles'',
arXiv:1103.1143, (2011).

\bibitem{B99} Brown 1999 ``Interlacing eigenvalues in time reversible Markov chains" 
\emph{Math. Op. Res.} {\bf 24}, 847 - 864(1999)

\bibitem{CMS} P. Collet, S.  Mart{\'\i}nez, J. San Mart{\'\i}n, ``Quasi-stationary distributions: Markov chains, diffusions and dynamical systems"
Springer Science \& Business Media 2012.

\bibitem{DS}
 J.N.Darroch, E.Seneta,
``On quasi-stationary distributions in absorbing discrete-time finite Markov chains'',
\emph{ J. Appl. Prob.} {\bf 2}, 88-100 (1965).

\bibitem{DF90}
P.Diaconis, J. A. Fill. 1990. ``Strong stationary times via a new form of duality", \emph{ Ann. Probab.}
{\bf 18}, no. 4, 1483?1522.

\bibitem{DM09}
P.Diaconis, L.Miclo, `` On Times to Quasi-Stationary for Birth and Death Processes",
\emph{ Journal of Theoretical Probability}, {\bf 22} (3) 558-586 (2009)

\bibitem{FMNS} R. Fern\'andez, F.\ Manzo, F.R.\ Nardi, E. Scoppola,
ï¿œAsymptotically exponential hitting times and metastability: a pathwise
ï¿œ approach without reversibilityï¿œ, 
Elettronic Journal of Probability 20 (2015) 122, 1-37

\bibitem{FMNSS} R. Fern\'andez, F.\ Manzo, F.R.\ Nardi, E. Scoppola, J. Sohier,
``Conditioned, quasi-stationary, restricted measures and escape from
metastable states'', \emph{Ann.Appl.Prob.}, {\bf 26} 760-793 (2016).

\bibitem{FL}
J.A.\ Fill, V.\ Lyzinski,
``Hitting times and interlacing eigenvalues: a stochastic appoach using intertwining'',
\emph{Journal of Theoretical Probability}, {bf 28},
Springer Science+Business Media New York 201210.1007/s10959-012-0457-9 (2012).


\bibitem{Kei79}
J.\ Keilson,
\emph{Markov Chain Models--Rarity and Exponentiality}, Springer-Verlag (1979).

\bibitem{LPW}
D.A.\ Levin, Y.\ Peres, E.L.\ Wilmer
\emph{ Markov Chains and Mixing Times}, AMS (2009).

\bibitem{OV}
E.\ Olivieri and M.E.\ Vares,
\emph{Large deviations and metastability}
Encyclopedia of Mathematics and its Applications, 100.
Cambridge University Press, Cambridge, (2005).


































%
%
%
%

\end{thebibliography}
\end{document}